\documentclass[12pt]{article}
\usepackage{amssymb,latexsym,mathrsfs}
\usepackage{graphicx}
\usepackage{fancyhdr}
\usepackage[colorlinks=true, linkcolor=blue, citecolor=blue,urlcolor=blue, breaklinks=true]{hyperref}
 \usepackage{lineno}
\usepackage{sectsty}
\usepackage{array}
\usepackage{booktabs}

\newtheorem{defs}{Definition}
\newtheorem{prop}{Proposition}
\newtheorem{theos}{Theorem}  

\begin{document}

\sectionfont{\large}
\subsectionfont{\normalsize}
\begin{center}
\textbf{\large Analysis and Applications of Delay Differential Equations in Biology and Medicine} 
\end{center}

\normalsize
\begin{center}
Majid Bani-Yaghoub,\\
Department of Mathematics and Statistics\\
                University of Missouri-Kansas City\\
                Kansas City, Missouri  64110-2499, USA\\
								baniyaghoubm@umkc.edu\\
\end{center}

\normalsize
\vspace{15pt}
\begin{center}
\noindent\textbf{Abstract}
\end{center}
The main purpose of this paper is to provide a summary of the fundamental methods for analyzing delay differential equations arising in biology and medicine. These methods are employed to illustrate the effects of time delay on the behavior of solutions, which include destabilization of steady states, periodic and oscillatory solutions, bifurcations, and stability switches. The biological interpretations of delay effects are briefly discussed.\\

\noindent \textbf{Math. Subj. Classification:} 37N25 (Dynamical systems in biology),  	37G15(Dynamical systems and ergodic theory )\\

\noindent \textbf{Key Words:} Allee Effect; Delay differential Equations; Stability Switch; Bifurcation\\

																%
					%
%
\normalsize
\section{Introduction}
\label{intro}
\label{sec:2}

The purpose of this paper is to describe some basic methods that are widely employed in the study of functional differential equations (FDEs) with a special interest on first and second order Delay Differential Equations (DDEs). Meanwhile, we mention some significant outcomes of the analysis of certain DDEs that may be used in the studies of biological systems.\\ 
\indent There are several great books  \cite{Hale 1993}, \cite{Antonevich}, \cite{Faria}, \cite{Erbe 1995}, \cite{Agarwal 2004}, \cite{Lebedev}, \cite{Skubach97}, \cite{Kolman99}, \cite{Green bk}, \cite{Kuang} in the fields of pure and applied mathematics devoted to the qualitative theory of differential equations with delays. Although some of them might seem very theoretical without giving an outline of the possible applications in biology, they are essential sources and references for work on FDEs. Namely, the second version of the book by J.Hale (coauthored by S.V. Lunel) \cite{Hale 1993} covers the basic theory of FDEs and also takes into account most of the fundamental achievements in the field. This volume refers in particular to the basic existence theory, properties of the solution map, Liapunov stability theory, stability and boundedness in general linear systems, behaviour near equilibrium and periodic orbits for autonomous retarded  equations, global properties of delay equations and FDEs on manifolds, which makes a great reference for studies related to FDEs.\\
\indent Since these basic developments of FDEs in 1970s, various approaches have been applied in the study of FDEs. For instance, the fundamental principles underlying the interrelations between c*-algebra and functional differential objects have been revealed in the book by A. Antonevich \cite{Antonevich}, where solvability conditions of various FDEs are investigated. In addition, the properties of solutions of FDEs have been examined with respect to oscillation theory in several study cases \cite{Green bk},  \cite{Erbe 1995}, \cite{Agarwal 2004}, \cite{Agarwal 2003} in which oscillatory and nonoscillatory properties of first, second and higher-order delay and neutral delay differential equations are addressed.\\
\indent While there are several contributions to the study of existence of solutions and also solvability of FDEs, the foremost concern of many applied mathematicians is the behaviour of the existing solutions of FDEs. Along with the new methods invented in the study of FDEs, some major tools such as method of characteristics or method of Liapunov functionals employed in global and local analysis of ODEs and PDEs have been extended to the analysis of FDEs. The present work covers some of the methods that have been discussed in the book by Gopalsamy \cite{Green bk} and the book by Kuang \cite{Kuang}.\\
%
%
%
\indent The classical stability theory of ODEs was generalized in the 1970s to investigate the stability of solutions of retarded functional differential equations (RFDEs) of the form
\begin{equation}
\label{eqCh3:genRFDE1}
\dot{x}(t)=f(t,x_{t}),
\end{equation}
where $x\in C([-\sigma-\tau,\sigma+A],\mathbb{R}^{n})$ with $\sigma\in\mathbb{R}$ $\tau,A\geq 0)$ and $t\in[\sigma,\sigma+A]$, $x_{t}\in C$ is defined as $x_{t}(\theta)=x(t+\theta)$ for $\theta\in[-\tau,0]$, $f:\mathbb{R}\times C([-\tau,0],\mathbb{R}^{n})\rightarrow\mathbb{R}^{n}$ .\\ Moreover, suppose that $f$ is uniformly continuous; then the stability of trivial solutions of system (\ref{eqCh3:genRFDE1}) is defined as follows.\\
\begin{defs}
Suppose $f(t,0)=0$ for all $t\in\mathbb{R}$. The solution $x=0$ of equation (\ref{eqCh3:genRFDE1}) is said to be
\begin{itemize}
	\item (i) stable if for any $\sigma\in \mathbb{R}, \epsilon>0,$ there is a $\delta=\delta (\epsilon,\sigma)$ such that $\phi\in \beta(0,\epsilon)$ implies $x_{t}(\sigma,\phi)\in\beta(0,\epsilon)$ for $t\geq\sigma$, where $\beta(0,\epsilon)$ is an open ball centered at the origin with radius $\epsilon$.\\
	\item (ii) asymptotically stable if it is stable and there is a $b_{o}=b_{o}(\sigma)>0$ such that $\phi\in\beta(0,b_{o})$ implies $x(\sigma,\phi)(t)\rightarrow0$ as $t\rightarrow\infty$.\\
	\item (iii) uniformly stable if the number $\delta$ in the definition is independent of $\sigma$.\\
	\item (iv) uniformly asymptotically stable if it is uniformly stable and there is a $b_{o}>0$ such that, for every $\eta>0$, there is a $t_{o}(\eta)$ such that $\phi\in\beta(0,b_{o})$ implies $x_{t}(\sigma,\phi)\in\beta(0,\eta)$ for $t\geq\sigma+t_{o}(\eta)$ for every $\sigma\in\mathbb{R}$.\\
\end{itemize}
\end{defs}
For a general solution $y(t)$ of system (\ref{eqCh3:genRFDE1}), the above-mentioned stability concepts are defined through the stability for the solution $z=0$ of the system
\begin{equation}
\label{eqCh3:genRFDE2}
\dot{z}(t)=f(t,z_{t}+y_{t})-f(t,y_{t}).
\end{equation}
\indent A usual manner in the study of FDEs is to investigate local stability analysis of some special solutions (e.g. trivial or constant solutions). For this purpose, the standard approach is to analyze the stability of equations linearized about the special solution. Hence, the stability of the special solution depends on the location of the roots of the related characteristic equation.
Despite the fact that RFDEs share many properties with ODEs (and also PDEs), we should emphasize that there are fundamental distinctions between the two theories. For instance, the linearized autonomous RFDEs define strongly continuous semi-groups on the phase space that are not analytic. In fact, the spectra of their generators consist of isolated eigenvalues with finite multiplicities. To explain this better, consider the linear delay differential equation

\begin{equation}
\label{eqCh2:linDDE1}
\dot{x}(t)=-\alpha x(t-\tau),
\end{equation}
which has a discrete delay term $\tau$. For simplicity let $\tau =1$. Then equation (\ref{eqCh2:linDDE1}) has the solution $x:t\rightarrow e^{\lambda t}$ if and only if the eigenvalue $\lambda$ is satisfied in the corresponding characteristic equation
\begin{equation}
\label{eqCh2:linDDE2}
\lambda+\alpha e^{-\lambda}=0.
\end{equation}
Nevertheless, only a finite number of eigenvalues may have a non-negative real part. Therefore, the center and unstable manifolds of the trivial solution are finite dimensional and the strongly continuous semi-group $\Gamma(t)=\exp(\alpha t)$ related to (\ref {eqCh2:linDDE1}) is not analytic.\\

\section{Method of Steps to Obtain Numerical Solutions}
Several methods have been proposed to solve systems of DDEs. Typical methods for solving DDEs are the method of characteristics, Laplace transforms and method of steps. Moreover, many DDE solvers have been developed since the early 1970s that use Runge-Kutta methods, Hermite interpolation and multistep methods for solving systems of DDEs. For instance, Matlab solver ``dde23''  is based on a third-order Runge-Kutta method that uses Hermite interpolation of the old and new solution and derivative to obtain an accurate interpolation. \\
Gathering all the methods for solving DDEs numerically or analytically is beyond the scope of the present work. In this section, we provide the method of steps that is commonly used for solving DDEs subject to an initial history function.\\
Let $s>0$, and $P(t)$ be a known function in $C([-\delta,0],\mathbb{R})$. Then the problem is to find a function $x(t), t\geq 0$ such that
\begin{equation}
\label{eq:StepMeth1}
\dot{x}(t)=f(t,x(t),x(t-\delta)),
\end{equation}
subject to the initial condition $x(t)=P(t)$ on $[-\delta,0]$.\\
\begin{itemize}
	\item Step 1: If $t\in[-\delta,0]$ then $x(t)=P(t)=x_{o}(t)$.
	\item Step 2: If $t\in[0,\delta]$ then $x(t-\delta)=x_{o}(t-\delta)$.\\ So we solve $\dot{x}(t)=f(t,x(t),x_{o}(t-\delta))$ which gives us the solution $x_{1}(t)$.
	\item Step 3: If $t\in[\delta,2\delta]$ then $x(t-\delta)=x_{1}(t-\delta)$.\\ So we solve $\dot{x}(t)=f(t,x(t),x_{1}(t-\delta))$ which gives us the solution $x_{2}(t)$.
\end{itemize}
Hence, for each interval we find a solution for equation (\ref {eq:StepMeth1}) and the general solution $x(t)$ includes all solutions $x_{o},x_{1},...$ defined in specific intervals.\\
For instance, let
\begin{equation}
\label{eq:StepMeth2}
\dot{y}(t)=y(t-1)-y(t), \mbox{ for }t>0,
\end{equation}
and for $t\in[-1,0]$,
\begin{equation}
\label{eq:StepMeth3}
y(t)=(t-1)^{2},
\end{equation}
(i.e. $\delta=1$ and initial history function $P(t)$ is given in (\ref{eq:StepMeth3})).\\
Step 1: If $t\in[-1,0]$, then $y(t)=(t-1)^{2}$.\\
\\
Step 2: If $t\in[0,1]$, then $t-1\in[-1,0]$, using equation (\ref{eq:StepMeth3}) we have
$$y(t-1)=(t-2)^{2}.$$
Hence the differential equation on $[0,1]$ is
$$\dot{y}(t)=(t-2)^{2}-y(t),$$
which has the solution
$$y_{1}(t)=\frac{1}{3}(t-2)^{3}+Ce^{-t}\mbox{ on }[0,1].$$\\
\\
Step 3: For $t\in[1,2], t-1\in[0,1]$, we solve the differential equation $\dot{y}(t)=y_{1}(t-1)-y(t)$ and this process continues until the desired time interval is approached.\\

\section{Analytical Methods to Study Delay Models}
\subsection{Method of Reduction to ODEs}
In some cases DDEs can be equivalent to systems of ODEs due to the special nature of the kernel functions in the integral terms. The so-called method of the ``chain trick'' was first introduced by D.M. Fargue \cite{Fargue} in 1973 and has been broadly used since (see for example \cite{Woerz-Busekros}, \cite{Macdonald1989}, \cite{Post and Travis}). We explain this method for the Lotka-Volterra sytem with distributed delays which is given by 
\begin{equation}
\label{eq:VolIntDiff1}
\frac{dx_{i}(t)}{dt}=x_{i}(t)\left(b_{i}+\sum^{n}_{j=1}a_{ij}x_{j}(t)+\sum^{n}_{j=1}b_{ij}\int^{t}_{-\infty}f_{ij}(t-s)x_{j}(s) ds\right),
\end{equation}
$$i=1,2,...,n$$
where $b_{i}, a_{ij}, b_{ij}(i,j=1,2,...,n)$ are real constants and $f_{ij}: [0,\infty)\mapsto [0,\infty)$ are continuous scalar functions known as delay kernels and normalized such that
\begin{equation}
\label{eq:VolIntDiff2}
\int^{\infty}_{0}f_{ij}(s)ds=1; i,j=1,2,...,n.
\end{equation}
With specific delay kernels $f_{ij}$, sufficient conditions for global asymptotic stability of system (\ref{eq:VolIntDiff1})-(\ref{eq:VolIntDiff2}) has been studied by A. Woerz-Busekros \cite{Woerz-Busekros}. By choosing the kernel functions of the form $f_{ij}(t)=\alpha e^{-\alpha t}, \alpha>0$, system (\ref{eq:VolIntDiff1}) is written as
\begin{equation}
\label{eq:VolIntDiff3}
\frac{dx_{i}(t)}{dt}=x_{i}(t)\left(b_{i}+\sum^{n}_{j=1}a_{ij}x_{j}(t)+\sum^{n}_{j=1}\beta_{ij}\alpha\int^{t}_{-\infty}e^{-\alpha(t-s)}x_{j}(s)ds\right),
\end{equation}
$$i=1,2,...,n;t>0,$$
where $b_{i}, a_{ij}, \beta_{ij}(i=1,2,...,n)$ are real constants and $\alpha$ is a positive constant.\\
Define a new set of variables $x_{n+j},j=1,2,...,n$ so that
\begin{equation}
\label{eq:VolIntDiff4}
x_{n+j}(t)=\alpha\int^{t}_{-\infty}e^{-\alpha(t-s)}x_{j}(s)ds; t>0.
\end{equation}
Using the product rule and the fundamental theorem of calculus we get that
\begin{equation}
\label{eq:VolIntDiff5}
\frac{dx_{n+j}(t)}{dt}=\alpha\left\{{x}_{j}(t)-x_{n+j}(t)\right\};  j=1,2,...,n.
\end{equation}
Thus, the system (\ref{eq:VolIntDiff3}) of $n$-integrodifferential equations becomes a system of $2n$ autonomous ordinary differential equations
\begin{equation}
\label{eq:autonODE1}
\frac{dx_{i}(t)}{dt}=x_{i}(t)\left(b_{i}+\sum^{n}_{j=1}a_{ij}x_{j}(t)+\sum^{n}_{j=1}\beta_{ij}x_{n+j}\right); i=1,2,...,n;
\end{equation}
$$\frac{dx_{n+j}(t)}{dt}=\alpha\left\{x_{j}(t)-x_{n+j}(t)\right\}; j=1,2,...n,.$$
If $x^{*}=(x^{*}_{1},x^{*}_{2},...,x^{*}_{n}),x^{*}_{1}>0, i=1,2,...,n$ is a solution of
\begin{equation}
\label{eq:autonODE2}
\sum^{n}_{j=1}(a_{ij}+\beta_{ij})x^{*}_{j}+b_{i}=0; i=1,2,...,n,
\end{equation}
then $(x^{*}_{1},x^{*}_{2},...,x^{*}_{n},x^{*}_{n+1},...,x^{*}_{2n}), x^{*}_{n+j}=x^{*}_{j}, j=1,2,...,n$ is a componentwise positive steady state of (\ref{eq:autonODE1}). Asymptotic stability of $(x^{*}_{1},...,x^{*}_{2n})$ for the system of ODEs (\ref{eq:autonODE1}) is equivalent to that of $(x^{*}_{1},...,x^{*}_{n})$ for DDEs system (\ref{eq:VolIntDiff3}).\\

\subsection{Method of Characteristics}
\label{Ch2 Loc Anal}
Local stability of a steady state solution of 
an ODE or PDE system is determined by linearizing the system at that steady state. The powerful Routh-Hurwitz criterion can be applied to the corresponding characteristic equations to determine if the real part of the roots are negative and if the steady state is stable. The method of characteristics has been extended to analyze the stability of DDEs. However, there are difficulties in applying such an extension. As manifested in the following example, in the presence of delay, the roots of the characteristic equation are functions of delays and hence, it is often a difficult task to apply the method of characteristics and the Routh-Hurwitz criterion to determine the local stability of steady state solutions. Consider the following prey-predator system with mutually interfering predators;
$$\frac{dx(t)}{dt}=x(t)\left[\gamma(1-\frac{x(t)}{k})-ay^{m}(t)\right],$$
\begin{equation}
\label{eqCh2:preyPred2}
\frac{dy(t)}{dt}=bx(t-\tau)y^{m}(t-\tau))-cy(t),
\end{equation}
where $a$, $b$, $c$ and $k$ are positive constants and $\tau\geq0$ is the discrete delay term, while $0<m<1$; $x(t)$ and $y(t)$ respectively denote the biomass of prey and predator populations. System (\ref{eqCh2:preyPred2}) has a positive steady state $E^{*}:=(x^{*},y^{*})$ satisfying
$$\gamma\left(1-\frac{x^{*}}{k}\right)=ay^{*m},$$
\begin{equation}
\label{eqCh2:preyPred4}
bx^{*}(y^{*})^{m-1}=c.
\end{equation}
Then by letting $x(t)=x^{*}+X(t)$ and $y(t)=y^{*}+Y(t)$ and linearizing (\ref{eqCh2:preyPred2}) around $E^{*}$, we arrive at
\begin{equation}
\label{eqCh2:linxpreyPred4}
\frac{dX(t)}{dt}= -\frac{\gamma}{k}x^{*}X(t)-amx^{*}(y^{*})^{m-1}Y(t),  \nonumber 
\end{equation}

\begin{equation}
\label{eqCh2:linypreyPred4}
\frac{dY(t)}{dt}=b(y^{*})^{m}X(t- \tau)+b(y^{*})^{m-1}x^{*}Y(t-\tau)-cY(t),
\end{equation}
which has the characteristic equation given by
\begin{equation}
\label{ch3:Dmatrix}
D(\lambda,\tau) =\bordermatrix{ & & \cr 
 & \lambda+ \frac{\gamma}{k}x^{*}& amx^{*}(y^{*})^{m-1} \cr 
 & -b(y^{*})^{m}e^{-\lambda \tau}&\lambda +c-bm(y^{*})^{m-1}x^{*}e^{-\lambda \tau}\cr
 }.                       
\end{equation}
\\  
  When $\tau=0$, $D(\lambda,0)$ is the usual quadratic equation of the form 
\begin{equation}
\label{eqCh2:preyPred8}
D(\lambda,0)=\lambda^{2}+p\lambda+q=0,
\end{equation}
where $p=c+\frac{\gamma}{k}x^{*}-bm(y^{*})^{m-1}x^{*}$  and\\ $q=\frac{\gamma}{k}x^{*}c-\frac{\gamma}{k}x^{*}m(y^{*})^{m-1}x^{*}+amx^{*}(y^{*})^{m-1}b(y^{*})^{m}$.\\
Using (\ref{eqCh2:preyPred4}), we can see that $p$ and $q$ are positive, which implies that (\ref{eqCh2:preyPred8}) has roots with negative real parts. Then, by the Routh-Hurwitz criterion, 
$E^{*}$ is locally asymptotically stable.\\
\indent When $\tau>0$, due to presence of terms with $e^{-\tau\lambda}$, the characteristic equation $D(\lambda,\tau)$ cannot be explicitly determined and hence, the linear stability analysis of the delay system (\ref{eqCh2:preyPred2}) via method of characteristics remains vague. In the best case for instance, one can establish sufficient conditions for the nonexistence of delay induced instability (i.e. conditions that the system (\ref{eqCh2:preyPred2}) remain stable near the steady state $E^{*}$ after inducing the delay $\tau>0$).\\
\indent However, one may bypass such difficulty by using the method of Liapunov functionals to obtain sufficient conditions for stability and instability of steady states of DDEs. Moreover, the stability results obtained in this way are often global.\\

\subsection{Method of Liapunov Functionals}
Consider the general retarded delay differential equation (RDDE),
\begin{equation}
\label{eqCh2:GenRDDE1}
\dot{x}(t)=f(t,x_{t}),
\end{equation}
where $x\in C([-\sigma-\tau,\sigma+A],\mathbb{R}^{n})$ with $\sigma\in\mathbb{R}$; $\tau,A\geq0$ and $t\in[\sigma,\sigma+A]$, $x_{t}\in C$ defined as $x_{t}(\theta)=x(t+\theta)$ for $\theta\in[-\tau,0]$, $f:\mathbb{R}\times C([-\tau,0],\mathbb{R}^{n})\rightarrow\mathbb{R}^{n}$ is uniformly continuous and $f(t,0)=0$.\\
Let $V:\mathbb{R}\times C\rightarrow\mathbb{R}$ be a continuous functional and $x(\sigma,\phi)$ be a solution of (\ref{eqCh2:GenRDDE1}) with initial value $\phi$ at $\sigma$ (i.e. there is an $A>0$ such that $x(\sigma,\phi)$ is a solution of (\ref{eqCh2:GenRDDE1}) on $[\sigma-\tau,\sigma+A)$ and $x_{\sigma}(\sigma,\phi)=\phi$).\\
Denote 
\begin{equation}
\label{eqCh2:GenRDDE2}
\dot{v}=V(t,\phi)=\overline{\lim}_{h\rightarrow0^{+}} \frac{1}{h}[V(t+h,x_{t+h}(t,\phi))-V(t,\phi)].
\end{equation}
The following theorem contains uniform (asymptotic) stability and boundedness results for the trivial solution of (\ref{eqCh2:GenRDDE1}).  \\

\begin{theos} 
\label{th1}
Let $u(s)$, $v(s)$, $w(s):\mathbb{R}^{+}\rightarrow\mathbb{R}^{+}$ be continuous and nondecreasing; $u(s)>0$, $v(s)>0$ for $s>0$ and $u(0)=v(0)=w(0)=0$.\\
The following statements are true:\\
(i) if there is a $V:\mathbb{R}\times C\rightarrow\mathbb{R}$ such that\\
$$u(\left|\phi(0)\right|)\leq V(t,\phi)\leq v(\left|\phi\right|),$$
$$\dot{v}(t,\phi)\leq-w(\left|\phi(0)\right|),$$
then $x=0$ (i.e. the trivial solution of (\ref{eqCh2:GenRDDE1})) is uniformly stable.\\
(ii) if in addition to (i) $\lim_{s\rightarrow+\infty} u(s)=+\infty$, then the solutions of (\ref{eqCh2:GenRDDE1}) are uniformly bounded (that is for any $\alpha>0$ there is a $\beta=\beta(\alpha)>0$ such that for all $\sigma\in \mathbb{R}$, $\phi\in C$ $\left\|\phi\right\|\leq\alpha$, we have $\left|x(\sigma,\phi)(t)\right|\leq\beta$ for all $t\geq\sigma$).\\
(iii) if in addition to (i), $w(s)>0$ for $s>0$, then $x=0$ is uniformly asymptotically stable.
\end{theos}
Therefore, the method includes the search of functionals $V$ satisfying the conditions of Theorem \ref{th1} to obtain stability for the trivial solution. For instance, the generalized form of Lotka-Volterra system (\ref{eq:VolIntDiff1})-(\ref{eq:VolIntDiff2}) is in the following form
\begin{equation}
\label{eq:GenVID1}
\frac{dx_{i}(t)}{dt}=x_{i}(t)\left(b_{i}+\sum^{n}_{j=1}a_{ij}x_{j}(t)+\sum^{n}_{j=1}b_{ij}x_{j}(t-\tau_{ij})+\sum^{n}_{j=1} c_{ij}\int^{t}_{-\infty} k_{ij}(t-s)x_{j}(s)ds\right),
\end{equation}
$$t>0; i=1,2,...,n;$$
with initial conditions 
\begin{equation}
\label{eq:GenVID2}
x_{i}(s)=\varphi_{i}(s)\geq0; s\in (-\infty,0); \varphi_{i}(0)>0; \sup_{s\leq0}\left|\varphi_{i}(s)\right|<\infty.
\end{equation}
Consider the Liapunov functionals $v(t)=v(t,x_{1}(\cdot),...,x_{n}(\cdot))$ defined by
\begin{eqnarray}
\label{eq:GenVID3}
v(t)&=&\sum^{n}_{i=1}\left|\log\left\{x_{i}(t)/x^{*}_{i}\right\}\right|+\sum^{n}_{j=1}\left|b_{ij}\right|\int^{t}_{t-\tau_{ij}}\left|x_{j}(s)-x^{*}_{j}\right|ds \\\nonumber
& &  +\sum^{n}_{j=1}\left|c_{ij}\right|\int^{\infty}_{0}\left|k_{ij}(s)\right|\left(\int^{t}_{t-s}\left|x_{j}(u)-x^{*}_{j}
\right|du\right)ds,  \mbox{ for }t\geq0.\\\nonumber
\end{eqnarray}
Then using Theorem \ref{th1} with a few sufficient conditions, it can be shown that all solutions of (\ref{eq:GenVID1}) subject to initial conditions (\ref{eq:GenVID2}) satisfy $\lim_{t\rightarrow\infty}x_{i}(t)=x^{*}_{i}; i=1,2.,...,n$, where $x^{*}=(x^{*}_{1},...,x^{*}_{n})$ is the positive system (\ref{eq:GenVID1}). Hence, $x^{*}$ is a global attractor and using the above method provides global asymptotic stability of $x^{*}$.


As stated below, the Liapunov functionals can also give sufficient conditions for the instability of the solution $x=0$ of a RDDE (\cite{Kuang} chapter 2):\\
\begin{theos} 
\label{th2}
Suppose $V(\phi)$ is a completely continuous scalar functional on $C$ and there exists a $\gamma >0$ and an open set $U$ in $C$ such that \\
(i) $V(\phi)>0$ on $U$, $V(\phi)= 0$ on the boundary of $U$; $0 \in cl( U \bigcap B(0,\gamma)$; \\
(ii) $V(\phi) \leq u(|\phi(0)|)$ on $U \bigcap B(0,\gamma)$;\\
(iii) $\underline{\dot{V}}(\phi) \geq w(|\phi(0)|)$ for  $ (t,\phi)\in [0,\infty) \times U \bigcap B(0,\gamma)$, where \\ 
$\underline{\dot{V}}(\phi) \equiv \underline{lim}\frac{1}{h}[V(x_{t+h}(t,\phi))-V(\phi)]$ as $ h \rightarrow 0^{+}$;\\
and where $u(s)$,$w(s)$ are continuous, positive and increasing for $s>0$, $u(0)=w(0)=0$. Then the trivial solution of (\ref{eqCh2:GenRDDE1}) is unstable.
\end{theos} 
The book by K. Gopalsamy \cite{Green bk} is a collection of different theorems and propositions on the stability of DDEs, many of them use the method of Liapunov functionals to establish a condition of local and global stability.\\
\indent In contrast to the favorable outcomes of the Liapunov functionals method mentioned above, there is a downside in employing this method to real problems arising from mathematical models. It is frequently quite demanding to find a Liapunov functional $V$ satisfying the conditions mentioned in Theorem \ref{th1} or \ref{th2}. Similarly, the method of DDEs stability analysis by employing Razumikin-type Theorems \cite{Kuang} suffers from difficulty of finding continuously differentiable functions that satisfy the conditions of the theorem. That is the reason many authors resort to the method of characteristic to obtain stability conditions for linear (or linearized) differential equations with discrete or distributed delays. However, in general, determining which of the methods is most advantageous over the others depends on the nature of the problem.\\

\subsection{Method of Hopf Bifurcation}
The classic Hopf bifurcation theory has been extended to systems of DDEs and also delayed PDEs by a number of authors. To explain this better, let us consider the linear system,
\begin{equation}
\label{eq1:LinearSys}
\frac{dx(t)}{dt}+ax(t)+bx(t-\tau)=0,
\end{equation}
where $b>a>0$. Then it has the characteristic equation
\begin{equation}
\label{eq2:LinearSys} 
\lambda+a+be^{-\lambda\tau}=0.
\end{equation}
If $\lambda=\mu+iw$ is a root of (\ref{eq2:LinearSys}), then so is $\bar{\lambda}=\mu-iw$. A Hopf bifurcation is subject to existence of a pair of pure imaginary eigenvalues in the case that $\mu=0$. Then, substituting $\lambda=iw$ into (\ref{eq2:LinearSys}) and solving for $\tau$ and $w$, we get that $w_{o}=\sqrt{b^{2}-a^{2}}$ and $\tau_{o}=\cos^{-1}(-a/b)/w_{o}$, where $\tau_{o}$ represents the Hopf bifurcation value and we have a case that is called ``delay induced bifurcation''. In particular, for values of $\tau$ near $\tau_{o}$, the trivial solution of (\ref{eq1:LinearSys}) is asymptotically stable for $\tau<\tau_{o}$ and it loses its stability when $\tau>\tau_{o}$. We have a similar situation for all $\tau_{o}+2k\pi, k=1,2,...$. Thus, for $\tau=\tau_{o}$ the linear variational system (\ref{eq1:LinearSys}) has periodic solutions with a period of  $\frac{2\pi}{w_{o}}$.\\
\indent A local bifurcation analysis can be conducted by perturbation methods (i.e. let $\tau=\tau_{o}+\epsilon$ with $0<\epsilon\ll1$ and substitute it in the characteristic equation) to demonstrate that small perturbations to the bifurcation value $\tau_{o}$ may destabilize the periodic solutions of the linear system (\ref {eq1:LinearSys}). The period of the exponentially growing unstable solutions of (\ref {eq1:LinearSys}) can be determined and the existence of nonlinear solutions near the perturbed bifurcation value can be established through the procedure of ``two-time asymptotic'' (see \cite{Murray 1984} for more details). The article by Mackey and Milton \cite{Murray I:343} and also \cite{Murray I:173} provides a good review of such analysis applied to the study of periodic dynamic diseases. Namely, Cheyne-Stokes respiration (i.e. human respiratory ailment) can be manifested by an alteration in the regular breathing pattern that can be presented by a nonlinear delay differential equation. Linearizing such equation around its steady state $x_{o}$ gives rise to the equation (\ref {eq1:LinearSys}), where $x(t)$ and $\tau$ are respectively the level of arterial carbon dioxide $CO_{2}$ and the time lag between the oxygenation of the blood in the lungs and monitoring by the chemoreceptors in the brainstem \cite{Murray I:173}. Ventilation of $CO_{2}$ in blood is related to $x(t)$ through the Hill function (\cite{Murray I}  Chapter 1) where the coefficient $b$ can be written as a product of a constant and evaluated at $x_{o}$ (i.e. $b=\beta v^{'}_{o}$). In order to be biologically meaningful, instead of investigating Hopf bifurcation due to changes of delay $\tau$, consider $v^{'}_{o}=\alpha$ to be a Hopf bifurcation value. Then it can be shown that small increases to the value $\alpha$ destabilizes the trivial solution of (\ref {eq1:LinearSys}) and results in an unstable steady state and a stable limit cycle with an approximate period of $4\tau$ (see \cite{Murray I} section 1.4) for more details). Therefore the period and volume of breathing may dramatically change if the gradient of the ventilation $v^{'}_{o}$ becomes too large.\\
Consider now the system,
\begin{equation}
\label{eq1:PertLinSys}
\frac{dx(t)}{dt}+ax(t)+bx(t-\tau)=f(x(t),x(t-\tau)),
\end{equation}
as a perturbation of (\ref{eq1:LinearSys}) (i.e. $f$ takes small values) and $\tau$ as a perturbation of $\tau_{o}$. Then similar to the classical Hopf Bifurcation Theory \cite{Hopf 1942}, \cite{Green bk}, the question arises of whether periodic solutions of system (\ref{eq1:LinearSys}) are stable under such perturbations or not. The perturbed equation (\ref{eq1:PertLinSys}) has a periodic solution with a period which is a perturbation of that of the linear approximation (\ref {eq1:LinearSys}). This has been investigated in several population dynamic models  \cite{Cushing1979},  \cite{Hale 1993}, \cite{Chow and Mallet-Peret 1977}, \cite{Stech 1979}, where in most cases, certain conditions on the parameter values are required to preserve stability of bifurcating periodic solutions under perturbations induced by delay. Chapter 2 of Gopalsamy's book \cite{Green bk} provides a self-contained demonstration of delay induced bifurcation to periodities of this type.\\

\subsection{Oscillatory and Nonoscillatory Methods}
\label{Ch2 Global Stab}
Oscillatory solutions of differential equations with or without delay have been frequently encountered in many biological processes described by a mathematical model. Chapters 7 - 9 of the book by J.D. Murray \cite{Murray I} provide a thorough background regarding biological and physiological oscillators studied via systems of differential equations. An oscillatory solution is generally defined as follows.\\
\begin{defs}
 A nontrivial solution $y$ is said to be \textit{oscillatory} if it has arbitrary large zeros for $t\geq t_{o}$, that is, there exists a sequence of zeros {$t_{n}$} (i.e $y(t_{n})=0)$ of $y$ such that $\lim_{n\rightarrow\infty} t_{n}=\infty$. Otherwise $y$ is said to be nonoscillatory.
\end{defs}
For instance, the second order delay differential equation,
\begin{equation}
\label{eqCh2:OsciDDE1}
y^{''}(t)+\frac{1}{2}y^{'}(t)-\frac{1}{2}y(t-\pi)=0,\mbox{ for }t\geq0,
\end{equation}
has oscillatory solution $y(t)=1-\sin(t)$ (i.e. it has an infinite sequence of multiple zeros).\\
\indent Recent developments in the oscillation theory of DDEs are presented in the book by R.P. Agarwal et al. \cite{Agarwal 2003}. In connection with wave profile equations, the existence of nonoscillatory solutions of second order DDEs has been examined in Chapter 5 of this book (see also  \cite{baniW}). Here, let us begin with the delayed logistic equation (Hutchinson's equation):
\begin{equation}
\label{eqCh2:Hutch1}
\dot{x}(t)=\gamma x(t)[1-x(t-\tau)/k],
\end{equation}
\begin{equation}
\label{eqCh2:Hutch2}
x(t)=\phi(t)\mbox{ on }[-\tau,0], \nonumber
\end{equation}
where $\phi(t)\in C([-\tau,0],\mathbb{R})$ is the initial history function.\\
We may nondimensionalize equation (\ref{eqCh2:Hutch1}) and reduce the number of parameters. In particular, let
 $\bar{y}(\bar{t})=-1+x(t)/k$ and $t=\tau\bar{t}$, then (\ref{eqCh2:Hutch1}) can be rewritten as
\begin{equation}
\label{eqCh2:Hutch3}
\frac{d}{d\bar{t}}\bar{y}(\bar{t})=-\gamma\tau\bar{y}(\bar{t}-1)[1+\bar{y}(\bar{t})].
\end{equation}
By dropping the bars from $\bar{y}$ and $\bar{t}$ and denoting $\alpha=\gamma\tau$, we have
\begin{equation}
\label{eqCh2:Hutch4}
\dot{y}(t)=-\alpha y(t-1)[1+y(t)],
\end{equation}
with a new initial history function $\tilde{\phi}\in C([-1,0],\mathbb{R})$.\\
Integrating from (\ref{eqCh2:Hutch4}), it can be observed that
\begin{equation}
\label{eqCh2:Hutch5}
1+y(t)=(1+y(t_{o}))\exp\left\{-\alpha\int^{t-1}_{t_{o}-1}y(\xi)d\xi\right\},
\end{equation}
which implies $1+y(t)>0$ as long as $y(\xi)$ exists on $[-1,t-1]$. It can be demonstrated \cite{Wright 1955} that the solution $y(t)=y(\widetilde{\phi})(t)$ of (\ref{eqCh2:Hutch4}) with initial history function $\widetilde{\phi}$ is bounded and asymptotically tends to the steady state $y(t)\equiv0$ of (\ref{eqCh2:Hutch4}) (i.e. $\lim_{t\rightarrow+\infty} y(t)=0$) if $\alpha\leq\frac{3}{2}$, $\phi(\theta)\geq-1$ and $\phi(0)>-1$.\\
Therefore, the positive steady state $x(t)\equiv k$ of (\ref{eqCh2:Hutch1}) with initial function $\phi(t)$ is globally asymptotically stable for delay $\tau\leq\frac{3}{2\gamma}$.\\
This is done by considering two cases for $y(t)$.
If $y(t)$ is nonoscillatory, then $y(t)>0$ or $y(t)<0$ for some $t\geq t_{o}\geq0$. Assume first that $y(t)>0$ for $t\geq t_{o}$; then, $\dot{y}(t)<0$ for $t\geq t_{o}+1$ from (\ref{eqCh2:Hutch4}) (since $1+y(t)>0$). Hence, $y(t)$ is strictly decreasing for $t\geq t_{o}+1$. There is a $c\geq 0$ such that $\lim_{t\rightarrow+\infty} y(t)=c$. And we must have $\lim_{t\rightarrow+\infty}\dot{y}(t)=0=-\alpha c(1+c)$. Therefore, $c=0$. The same conclusion holds for $y(t)<0$ for $t\geq t_{o}$.\\
In addition, if $y(t)$ is oscillatory then the global stability will be derived by using basic calculus and local maximum minimum properties of $y(t)$.\\
The result can be improved to $\tau\leq\frac{37}{24\gamma}$ at the cost of considerable elaboration \cite{Kuang}. Nevertheless,  the attempt to show $\alpha<\frac{\pi}{2\gamma}$ was not successful \cite{Gopalsamy1986}.\\
Equation (\ref{eqCh2:Hutch1}) for single species growth can be generalized to the following first order scalar non-autonomous delay equation with negative feedbacks \cite{Haddock and Kuang 1992}.
\begin{equation}
\label{eqCh2:scalar1}
\dot{x}(t)=\int^{t}_{t-\tau(t)}\sum^{n}_{i=1} f_{i}(t,x(s))d\mu_{i}(t,s),
\end{equation}
where $\tau(t)>0$, $f_{i}(t,x)$ and $\tau(t)$ are continuous with respect to their arguments and $\mu_{i}(t,s)$ is continuous with respect to $t$, nondecreasing with respect to $s$ and is defined for all $(t,s)\in\mathbb{R}^{2}$. 
Then with a few modifications to the previous method, sufficient conditions for global stability of the trivial solution of (\ref {eqCh2:scalar1}) (with respect to appropriate initial functions) are established \cite{Haddock and Kuang 1992}.\\

\section{Applications in Biology and Medicine}
\subsection{Destabilizing Effect of Delay}
\indent In this section we discuss the effects of delay on the behavior of solutions for models that have been developed for study in different problems in biology and medicine. As mentioned in Subsection \ref{Ch2 Global Stab}, in the study of the delay effects on systems of differential equations, there are many articles that consider delay as a small perturbation to the system. Then the perturbation methods can be used to take advantage of already known results in non-delayed systems of ODEs or PDEs. For instance, it can be demonstrated that the small delays have no influence on the qualitative behaviour of the solution of the delayed logistic equation (\ref{eqCh2:Hutch1}); whereas, large delays destabilize its positive steady state. Such small delays with negligible effects on the behaviour of the solution are often referred to as harmless delay \cite{Gopalsamy 1983Red}. One may think that sufficiently small delays are always harmless and can be ignored in the model analysis but this is not so.\\
\indent A counter example may be found in the book by Hale (\cite{Hale 1993} version 1977, page 28)  where the trivial solution of
\begin{equation}
\label{eq1:HaleTrivSol}
\dot{x}(t)+2\dot{x}(t)=-x(t),
\end{equation}
is asymptotically stable, but the trivial solution of
\begin{equation}
\label{eq2:HaleTrivSol}
\dot{x}(t)+2\dot{x}(t-\tau)=-x(t),
\end{equation}
is unstable to any positive delay $\tau$. Other examples of this type may be found in Kolmanovskii and Nosov and the invariant systems studies by Shipanov \cite{Shipanov1939}. Moreover, the destabilizing effect of delay can be seen in general scalar neutral differential equations with a single delay $\tau\geq0$
\begin{equation}
\label{eqCh2:neutDiff1}
\sum^{n}_{k=0} a_{k}\frac{d^{k}}{dt^{k}}x(t)+\sum^{n}_{k=0}b_{k}\frac{d^{k}}{dt^{k}}x(t-\tau)=0.
\end{equation}
Then using the method of characteristic, it can be demonstrated (\cite{Kuang} chapter 3) that 
the trivial solution of equation (\ref{eqCh2:neutDiff1}) loses its stability for any $\tau>0$ when $\left|b_{n}\right|>0$. \\

\subsection{Oscillation or Nonoscillation Affected by Delay}
\indent The effect of delay on the oscillatory and nonoscillatory behaviour of delay differential equation
\begin{equation}
\label{eqCh2:OsciDDE2}
\frac{dx(t)}{dt}+ax(t-\tau)=0,
\end{equation}
can be seen in the proposition that follows.\\
\begin{prop}
Let $a\in(0,\infty)$ and $\tau\in(0,\infty)$. Then all nontrivial solutions of (\ref{eqCh2:OsciDDE2}) are oscillatory if
\begin{equation}
\label{eqCh2:OsciDDE3}
ae\tau>1,
\end{equation}
and (\ref{eqCh2:OsciDDE2}) has a nonoscillatory solution if
\begin{equation}
\label{eqCh2:OsciDDE4}
ae\tau\leq1.
\end{equation}
\end{prop}
Such a result is a very special case of a large class of DDEs studied by J. Yan \cite{Yan 1987}. Another example is oscillations in a Lotka-Volterra system that has been well investigated by Gopalsamy \cite{Gopalsamy1991}. So far, we have seen that delay may have an effect on the stability of steady states, asymptotic behaviour of trivial solutions and oscillatory (or nonoscillatory) behaviour of solutions. In the following we will observe that delay may cause phenomena called ``stability switches. ''\\

\subsection{Stability Switches}
The other phenomena to mention in this section are the stability switches due to changes of delay. Starting with a small delay, as the length of the delay increases, the trivial solution of DDEs can gain or lose its linear stability. Such phenomena are often called stability switches. There are plenty of studies providing sufficient conditions for the existence or nonexistence of stability switches (see \cite{Freedman and Gopalsamy 1988} and the references therein). Here, consider the system of DDEs
$$\frac{dx_{1}(t)}{dt}=r_{1}x_{1}(t)\left[\frac{k_{1}+\alpha_{1}x_{2}(t-\tau_{2})}{1+x_{2}(t-\tau_{2})}-x_{1}(t)\right],$$
\begin{equation}
\label{Cheq:linDDE2}
\frac{dx_{2}(t)}{dt}=r_{2}x_{2}(t)\left[\frac{k_{2}+\alpha_{2}x_{1}(t-\tau_{1})}{1+x_{1}(t-\tau_{1})}-x_{2}(t)\right],
\end{equation}
where $\tau_{1},\tau_{2}\geq0$ and $\tau_{1}+\tau_{2}>0$.\\
The system (\ref{Cheq:linDDE2}) indicates that the mutualistic or cooperative effects are not realized instantaneously but take place with time delays. Linearizing the system around its positive steady state $x^{*}$ and solving the corresponding characteristic equation for $\lambda=\alpha+i\beta$, the necessary and sufficient conditions can be found for nonexistence of stability switches. In particular we have the following theorem (\cite{Green bk}, Section 3.3).\\
\begin{theos}
Assume that $r_{i}, k_{i},\alpha_{i}>0$ and $\alpha_{i}>k_{i}$ for $i=1,2$; then the positive steady state $x^{*}$ of system (\ref{Cheq:linDDE2}) is linearly asymptotically stable absolutely in delays (i.e. delay induced stability switches cannot occur and $N^{*}$ is asymptotically stable for all delays).
\end{theos}

\subsection{Conditions for Delay Independent Stability}
Despite the fact that time delays are often thought to have destabilizing or stability switching effects, we may observe cases where local stability of a delay system is not affected by delay at all. In particular, the following theorem gives sufficient conditions for delay independent local stability of a steady state of a delay model. \\
Let $x_{1}(t)$ and $x_{2}(t)$ denote the population densities of two species competing for a common pool of resources in a temporally uniform environment; let $b_{i}$ and $m_{i}$ $(i=1,2)$ denote the respective density dependent birth and death rates (see \cite{Miller 1976}, \cite{Brian 1956}, \cite {Gopalsamy 1984b} for an extensive discussion of competition processes). Let $\tau_{ij}$ $(i,j=1,2)$ be a set of nonnegative constants with $\tau=\max\left\{\ \tau_{ij}|i,j=1,2\right\}$ so that the population densities are governed by
$$\frac{dx_{1}(t)}{dt}=b_{1}(x_{1}(t-\tau_{11}))-m_{1}(x_{1}(t),x_{2}(t-\tau_{12})),$$
\begin{equation}
\label{eqCh2:popDens2}
\frac{dx_{2}(t)}{dt}=b_{2}(x_{2}(t-\tau_{22}))-m_{2}(x_{1}(t-\tau_{21}),x_{2}(t)),
\end{equation}
with initial population size
$$x_{i}(s)=\phi_{i}(s)>0; s\in[-\tau,0]; i=1,2,$$
\begin{equation}
\label{eqCh2:popDens4}
\phi_{i}\in C([-\tau,0],\mathbb{R}^{+}),\phi_{i}\neq 0\mbox{ on }[-\tau,0], i=1,2.
\end{equation}
The following assumptions on the birth and death rates are made for the system of DDEs (\ref{eqCh2:popDens2}):\\
(i) $b_{i}$, $m_{i}(i=1,2)$ are continuous with continuous partial derivatives for all $x_{i}\geq0(i1,2)$; also we assume
\begin{equation}
\label{eqCh2:popDens5}
\frac{\partial b_{i}}{\partial x_{i}}>0; \frac{\partial m_{i}}{\partial x_{j}}>0\mbox{ for }x_{i}>0; i,j=1,2;
\end{equation}
(ii)
\begin{equation}
\label{eqCh2:popDens6}
b_{i}(0)=0; m_{i}(0,x_{2})\equiv0;
\end{equation}
(iii) for some $x^{*}_{1}>0$, $x^{*}_{2}>0$ we have
$$b_{1}(x^{*}_{1})-m_{1}(x^{*}_{1},0)=0,$$
\begin{equation}
\label{eqCh2:popDens7}
b_{2}(x^{*}_{2})-m_{2}(0,x^{*}_{2})=0;
\end{equation}
(iv) there exist positive constants $\delta_{1}$,$\delta_{2}$ such that
$$b_{1}(\delta_{1})-m_{1}(\delta_{1},x_{2})<0,$$
and
\begin{equation}
\label{eqCh2:popDens8}
b_{2}(\delta_{2})-m_{2}(x_{1},\delta_{2})<0,
\end{equation}
for $x_{1}\geq0$, $x_{2}\geq0$;\\
(v) for the positive steady state $(\alpha,\beta)$ of (\ref{eqCh2:popDens2}), we have
$$b_{1}(\alpha)-m_{1}(\alpha,\beta)=0,$$
$$b_{2}(\beta)-m_{2}(\alpha-\beta)=0;$$
(vi) $$\frac{\partial m_{1}}{\partial x_{1}}>\frac{\partial b_{1}}{\partial x_{1}}+\frac{\partial m_{2}}{\partial x_{1}},$$
$$\frac{\partial m_{2}}{\partial x_{2}}>\frac{\partial b_{2}}{\partial x_{2}}+\frac{\partial m_{1}}{\partial x_{2}}.$$
\begin{theos} 
Assume that the conditions (i)-(vi) hold for the two species competition model (\ref {eqCh2:popDens2}) with delays in production and interspecific competitive destruction, then the positive steady state $(\alpha,\beta)$ of (\ref{eqCh2:popDens2}) is (locally) asymptotically stable for all delays $\tau_{ij}\geq0; i,j=1,2$.
\end{theos}
\subsection{Stability Conditions for DDEs}
There are numerous articles employing the above-mentioned methods to establish the conditions for local or global stability of solutions of different delay models. Here, we provide two of many outcomes of this type available for different DDEs.\\
\indent Consider the following second order differential equation with finite number of discrete delays 
\begin{equation}
\label{eqCh2:linWave1}
x^{''}(t)+a_{1}x^{'}(t)+a_{o}x(t)=\sum^{n}_{j=1} b_{j}x(t-\tau_{j}),
\end{equation}
where $a_{o}$, $a_{1}$ and $\tau_{j}\geq0$; $b_{j}\in\mathbb{R}$ for $j=1,...,n$ and $a_{o}>\sum^{n}_{j=1} \left|b_{j}\right|$.\\
\indent Then the following theorem by G. Stepan \cite{Stepan1989} gives the conditions for uniform asymptotic stability of the trivial solution of (\ref {eqCh2:linWave1}).

\begin{theos} 
The trivial solution of (\ref {eqCh2:linWave1}) is uniformly asymptotically stable for all values of $\tau_{j}\geq0$ if either
\begin{equation}
\label{eqCh2:linWave2}
a_{1}>\frac{\sum^{n}_{j=1}\left|b_{j}\right|}{\left(a_{o}-\sum^{n}_{j=1}\left|b_{j}\right|\right)^{\frac{1}{2}}},
\end{equation}
or
\begin{equation}
\label{eqCh2:linWave3}
a_{1}>\sum^{n}_{j=1}\left|b_{j}\right|\tau_{j}.
\end{equation}
\end{theos}
In the case of distributed delay, the equation (\ref {eqCh2:linWave1}) is changed to 
\begin{equation}
\label{eqCh2:linWave4}
x^{''}(t)+a_{1}x^{'}(t)+a_{o}x(t)=\int^{0}_{-\tau} x(t+\theta) d\eta(\theta),
\end{equation}
where $\int^{0}_{-\tau} \left|d\eta(\theta)\right|=\eta<+\infty$, $a_{o}>\eta$ and there is a $v>0$ such that $\int^{0}_{-\tau} e^{-v\theta}\left|d\eta(\theta)\right|<+\infty$.\\
Then we have the following theorem:\\
\begin{theos} 
The trivial solution of (\ref{eqCh2:linWave4}) is uniformly asymptotically stable if either 
\begin{equation}
\label{eqCh2:linWave5}
a_{1}>\frac{\eta}{(a_{o}-\eta)^{\frac{1}{2}}},
\end{equation}
or
\begin{equation}
\label{eqCh2:linWave6}
a_{1}>\int^{0}_{-\tau}\left|\theta\right|d\eta(\theta).
\end{equation}
\end{theos}

\section{Discussion}
Linearization at a steady state is one of the main tools in studying continuous mathematical delay models representing population or epidemic dynamics \cite{Bani Miami}, \cite{Bani TH},  \cite{Bani08},  \cite{Bani II}, \cite{bani10}. For example, several theorems establish conditions for asymptotic stability of the trivial solution (i.e. zero solution) of a delay differential equation or a neutral functional differential equation (NFDE) through the analysis of their characteristic equations. Nevertheless, linearization results often provide information about the behavior of solutions only near a steady state. In general, questions such as existence and global stability of periodic orbits or oscillations in solutions of a continuous delay model can be answered by employing asymptotic methods, bifurcation analysis, method of Liapunov functionals and other known methods in global analysis of a system of FDEs or PDEs. The global stability analysis of steady states has been the focus of many researchers examining various delay differential population models. Namely, it is often desirable to obtain sufficient conditions for the global asymptotic stability of the positive steady state of a nonlinear differential equation. Furthermore, the global existence of periodic solutions and also chaotic behavior induced by delay have been investigated in quite a few studies \cite{Wright 1955}, \cite{Haddock and Kuang 1992}, \cite{Green bk}, \cite{Kuang},  \cite{Li and Yorke 1975}\\
\indent In the study of qualitative changes to DDEs and NFDEs due to changes of discrete (and also distributed) delays, several authors (e.g. \cite{cooke grossman 1982}, \cite{cooke Driesches}, \cite{FreedKuang 1991}) have encountered the stability switches that may take place for a trivial solution of a non-autonomous DDE. In particular, the stability of the trivial solution will be affected through increases of the delay length. Moreover, delay may induce destabilizing effects, oscillatory effects or no qualitative effect at all. \\
\indent Providing the two cases of oscillatory and nonoscillatory among the necessary transformations can be effective in global analysis of DDEs and detecting the sufficient conditions of global stability of trivial solutions. Moreover, the method of Liapunov functionals and also Razumikhin-type theorems have been frequently applied to the study of global behavior of solutions. For instance,  Razumikhin functions can be used to show that \cite{Gopalsamy Ladas 1990} all positive solutions of (2.30) are attracted by the steady state $x^{*}$  when $cx^{*}> b$ and $\tau < 1/r$ ($r$ is a constant). Using the above-mentioned methods, the global and local analysis of several mathematical models in biology have been investigated in a number of studies \cite{hadeler and tomiuk 1977}, \cite{mallet-paret and Nussbaum 1986}, \cite{gopalsamy et al 1989}, \cite{gopalsamy et al 1990b},  \cite{kuang and smith 1991b}, \cite{Tang and kuang 1992a}, \cite{leung 1979}, \cite{leung and zhou 1988}.

Outcomes and analysis of delay models \cite{bani4}, \cite{bani1} are biologically interpreted according to each research project in which they have been applied. In the following we provide two well known effects in population biology that have been widely studied through systems of DDEs.
\subsection{Allee effect} Global behaviors of the solutions can biologically be interpreted in distinct ways. The so-called Allee effect \cite{Allee 1927}, \cite{Allee 1933} relates to a population that has a maximal per capita growth rate at intermediate density. When the population becomes too large, the positive feedback effect of aggregation and cooperation may then be dominated by density dependent stabilizing negative feedback effect due to intraspecific competition arising from excessive crowding and the ensuing shortage of resources. \\
These processes have been studied \cite{bani2}, \cite{bani13}  through global (and local) analysis of several models such as the following Lotka-Volterra type single species population growth \cite{Gopalsamy Ladas 1990}
\begin{equation}
\label{eqCh2:scalar2}
\dot{x}(t)=x(t)[a+bx(t-\tau)-cx^{2}(t-\tau)],
\end{equation}
with $x(t)=\phi(t)\geq0$, $t\in[-\tau,0)$ and $\phi(0)>0$, where $a,c>0$, $b\in\mathbb{R}$ and $\phi\in C([-\tau,0],\mathbb{R})$. \\
When $\tau=0$ and $b>0$, the system exhibits the Allee effect. Moreover, equation (\ref{eqCh2:scalar2}) has a unique positive equilibrium $x^{*}=\frac{1}{2c}(b+\sqrt{b^{2}+4ac})$.\\
Then the transformation $x(t)=x^{*}[1+y(t)]$ reduces equation (\ref{eqCh2:scalar2}) to 
\begin{equation}
\label{eqCh2:scalar3}
\dot{y}(t)=-\alpha(t)y(t-\tau), t\geq0,
\end{equation}
where $\alpha(t)=[(2cx^{*}-b)x^{*}+c(x^{*})^{2}y(t-\tau)][1+y(t)]$.\\
Note that such a transformation is required to make the derivative negative. In this way it can be shown that the trivial solution is a global attractor of all nonoscillatory solutions of (\ref {eqCh2:scalar2}) when $2cx^{*}-b>0$. In fact let $y(t)>0$ be a nonoscillatory solution, then from (\ref{eqCh2:scalar3}) we get that $\dot{y}(t)<0$. Since $y(t)$ is nonoscillatory, it is concluded that $\lim_{t\rightarrow\infty} y(t)\geq0$ which if strictly greater than zero, by letting $2cx^{*}-b>0$, we have $\lim_{t\rightarrow\infty} \dot{y}(t)<0$. However, the last inequality implies that $\lim_{t\rightarrow\infty} y(t)=-\infty$ which is a contradiction. The case $y(t)<0$ is similar to this argument.\\
By using the local maximum or minimum properties the global upper and lower bounds for oscillatory solutions of (\ref {eqCh2:scalar2}) are in the form of
$$e^{-M\tau}\leq1+y(t)\leq e^{Lx^{*}\tau}\mbox{ for }t\geq T,$$
where $y(t)$ is an oscillatory solution; $L,M$ and $T$ are constants.\\

\subsection{Permanence} While the local and global qualitative properties of a given system of DDEs or delayed PDEs are crucial to mathematically analyze and predict different phenomena in population biology and epidemiology, an important fundamental property to consider is the permanence (persistence) of the system in the long run. In particular, the question is whether the involved populations and/or epidemics will remain permanently in coexistence or one of them will finally survive at the expense of the other's extinction. Permanence of Lotka-Volterra type systems with delays \cite{Cao and Gard 1992},\cite{Wang and Ma 1991}, \cite{Burton and Hutsen 1989}, \cite{Kuang and Tang 1992a}, permanence of delayed Kolmogorove-type systems \cite{Cao et al. 1992b}, \cite{Gopalsamy Ladas 1990}, and the uniform persistence of functional differential equations are the typical works on permanence of DDEs. Note that, the general persistence theory of ODEs (see for example \cite{Hale and Waltman 1989}) is a prerequisite to these studies.\\

In conclusion, the present work provides the summary of the basic tools that are used in the studies of delay differential equations. More sophisticated tools \cite{bani11} \cite{bani2}  \cite{bani12}  have been developed based on the characteristics of each specific model.


\begin{thebibliography}{999}



\bibitem{Allee 1927} W.C. Allee (1927) Animal aggregations. Quart. Rev. Biol. 2, 367--398.

\bibitem{Allee 1933} W.C. Allee (1933) Animal aggregations: A Study in General Sociology. Chicago Univ. Press: Chicago.







\bibitem{Antonevich} Antonevich, Anatolii. (1998) Functional Differential Equations.  CRC Press.

\bibitem{Agarwal 2003} Ravi P. Agarwal, Said R. Grace, Donal O'Regan. (2003) Oscillation Theory for Second Order Dynamic Equations.  CRC Press.


\bibitem{Agarwal 2004} Ravi P. Agarwal, Martin Bohner, Wan-Tong Li.(2004) Nonoscillation and Oscillation: Theory for Functional Differential Equations.  Marcel Dekker, Inc.: New York.














\bibitem{Bani Miami} M. Bani-Yaghoub and D.E. Amundsen (2006) Turing-type instabilities in a mathematical model of Notch and Retinoic Acid pathways. WSEAS Transactions on Biology and Biomedicine 3(2), 89--96.

\bibitem{Bani TH} M. Bani-Yaghoub (2006) A mathematical approach to axon formation in a network of signaling molecules for N2a cells. M.Sc. Thesis, Carleton University, Ottawa, Canada. 

\bibitem{Bani08} M. Bani-Yaghoub and D.E. Amundsen (2008) Study and Simulation of Reaction–Diffusion Systems Affected by Interacting Signaling Pathways Acta biotheoretica 56 (4), 315-328

\bibitem{Bani II} M. Bani-Yaghoub (2010) Wave Solutions of Nonlocal Delayed Reaction-diffusion Equations, PhD. Thesis, Carleton University, Ottawa, Canada. 

\bibitem{bani10} M. Bani-Yaghoub and D.E. Amundsen (2010) Dynamics of Notch Activity in a Model of Interacting Signaling Pathways. Bulletin of mathematical biology 72 (4), 780-804

\bibitem{bani4} M. Bani-Yaghoub, G. Yao (2014) Modeling and Numerical Simulations of Single Species Dispersal in Symmetrical Domains,  International Journal of Applied Mathematics, 27(6), 525--547. 

\bibitem{bani1}  M. Bani-Yaghoub, D.E. Amundsen (2015) Oscillatory traveling waves for a population diffusion model with two age classes and nonlocality induced by maturation delay. Comput. Appl. Math. 34(1), 309--324.

\bibitem{bani11} M. Bani-Yaghoub (2017) Approximating the traveling wavefront for a nonlocal delayed reaction-diffusion equation, Journal of Applied Mathematics and Computing,  53 (1), 77--94

\bibitem{baniW} M. Bani-Yaghoub (2017) Introduction to Delay Models and Their Wave Solutions. arXiv preprint 1-20

\bibitem{bani2} M. Bani-Yaghoub, G. Yao,  M. Fujiwara, D.E. Amundsen (2015) Understanding the interplay between density dependent birth function and maturation time delay using a reaction-diffusion population model,  Ecological Complexity, 21 , 14--26. 

\bibitem{bani12} M. Bani-Yaghoub (2016) Approximate Wave Solutions Of Delay Diffusive Models Using A Differential Transform Method Applied Mathematics E-Notes 16, 99-104 

\bibitem{bani13} M. Bani-Yaghoub, G. Yao, H. Voulov (2016) Existence and stability of stationary waves of a population model with strong allee effect, Journal of Computational and Applied Mathematics,  307, 385-393



\bibitem{Brian 1956} M.V. Brian (1956) Exploitation and interference in interspecies competition. J. Anim. Ecol. 25, 339--347.





\bibitem{Burton and Hutsen 1989} T.A. Burton and V. Hutson (1989) Repellers in systems with infinite delay. J. Math. Anal. Appl. 137, 240--263.


\bibitem{Cao et al. 1992b} Y. Cao, J.-P. Fan and T.C. Gard (1992b) Uniform persistence for population interaction models with time delay. Applicable Analysis 51, 197 - 210.

\bibitem{Cao and Gard 1992} Y. Cao and T.C. Gard (1993) Uniform persistence for population models with time delay using multiple Lyapunov functions. Differential and Integral Equations 6(4), 883--898. 








\bibitem{Chow and Mallet-Peret 1977} S.N. Chow and J. Mallet-Peret (1977) Integral averaging and Hopf's bifurcation. J. Diff. Eqns. 26, 112--159.



\bibitem{cooke Driesches} K.L. Cooke and P. van der Driesche (1986) On zeros of some transcendental equations. Funkcialaj Ekvacioj 29, 77--90.

\bibitem{cooke grossman 1982} K.L. Cooke and Z. Grossman (1982) Discrete delay, distributed delay and stability switches. J. Math. Anal. Appl. 86, 592--627.


\bibitem{Cushing1979} J.M. Cushing (1979) Volterra inetegrodifferential equations in population dynamics. In: Mathematics in Biology. Proc. C.I.M.E.: Italy.






\bibitem{Erbe 1995} Lynn H. Erbe, Bing-Gen Zhang, Qingkai Kong. (1995) Oscillation Theory for Functional Differential Equations.  CRC Press.


\bibitem{Faria} T. Faria, P. Freitas. (2001) Topics in Functional Differential and Difference Equations. 
American Mathematical Society: Providence, RI.



\bibitem{Fargue} Fargue, D. (1973). Réductibilité des systèmes héréditaires à des systèmes dynamiques. Compt. Rend. Acad. Sci., B277, 471–473 




\bibitem{Freedman and Gopalsamy 1988} H.I. Freedman and K. Gopalsamy (1988) Nonoccurrence of stability switching in systems with discrete delays. Canad. Math. Bull. 31, 52--58.

\bibitem{FreedKuang 1991} Freedman, H.I. and Kuang, Y. (1991) Stability switches in linear scalar neutral delay equations. Funkcialaj Ekvacioj. 34, 187--209.





\bibitem{Murray I:173} L. Glass and M.C. Mackey (1988) From Clocks to Chaos: The rhythm of Life. Princeton Univ. Press:Princeton, NJ.



\bibitem{Gopalsamy 1983Red} K. Gopalsamy (1983) Harmless delays in model systems. Bull. Math. Biol. 45, 295--309.

\bibitem{Gopalsamy 1984b} K. Gopalsamy (1984b) Delayed responses and stability in two-species systems. J. Austral. Math. Soc. Ser. B 25, 473--500.

\bibitem{Gopalsamy1986} K. Gopalsamy (1986) On the global attractivity in a generalized delay-logistic differential equation. Proc. Camb. Phil. Soc. 100, 183--192.

\bibitem{Gopalsamy1991} K. Gopalsamy (1991) Oscillations in Lotka Volterra systems with several delays. J. Math. Anal. Appl. 159, 440--448.

\bibitem{Green bk} K. Gopalsamy (1992) Stability and Oscillations in Delay Differential Equations of Population Dynamics. Kluwer Academic Publishers: Dordrecht, The Netherlands.

\bibitem{gopalsamy et al 1989} K. Gopalsamy, M.R.S. Kulenovic, G. Ladas, (1989) Oscillation and global attractivity in respiratory dynamics, Dynam. Stab. Systems 4  131-139. 

\bibitem{gopalsamy et al 1990b} Gopalsamy, MRS Kulenovic, G. Ladas,(1990) Oscillation and global attractivity in models of hematopoiesis, J. Dyn. Diff. Eqns. 2  117-132.

\bibitem{Gopalsamy Ladas 1990} K. Gopalsamy and G. Ladas (1990) On the oscillation and asymptotic behavior of $N'(t)=N(t)[a+bN(t-\tau)-cN^{2}(t-\tau)]$. Quart. Appl. Math. 3, 433--440. 




 










\bibitem{Haddock and Kuang 1992} J.R. Haddock and Y. Kuang (1992) Asymptotic theory for a class of nonautonomous delay differential equations. J. Math. Anal. Appl. 168, 147--162.


\bibitem{hadeler and tomiuk 1977} Hadeler, K. P., and Tomiuk, J. (1977). Periodic solutions of difference-differential equations. Arch. Rat. Mech. Anal. 65, 87-95. 


\bibitem{Hale 1993}  J.K. Hale and L.S.M. Verduyn (1993) Introduction to Functional Differential Equations, Springer.

\bibitem{Hale and Waltman 1989} J.K. Hale and P. Waltman (1989) Persistence in infifnite-dimensional systems. SIAM J. Math. Anal. 20, 388--395.






\bibitem{Hopf 1942} E. Hopf (1942) (Bifurcation of a periodic solution from a staionary solution of a system of differential equations) Translated from German, Ber. Math. Phys. Klasse de Sachs. Akad. Wiss. Leipzig 94, 3--22.















\bibitem{Kolman99} Vladimir Borisovich Kolmanovskii, Anatolii Dmitrievich Myshkis. (1999)Introduction to the Theory and Applications of Functional Differential Equations.   Springer.









\bibitem{Kuang} Yang Kuang. (1993) Delay Differential Equations with Applications in Population Dynamics.  Academic Press, Inc.: SAn Diego.




\bibitem{kuang and smith 1991b} Y. Kuang and H. L. Smith,(1992) Convergence in Lotka-Volterra-type delay systems without in- stantaneous feedbacks, Proc. Roy. Soc. Edinburgh, 123 Sect.


\bibitem{Kuang and Tang 1992a} Y. Kuang and B.R. Tang (1993) Uniform persistence in nonautonomous delay differential Kolmogorov type population models. Rocky Mountain J. Math. 


\bibitem{Tang and kuang 1992a} Y. Kuang and B. R. Tang (1994)  Uniform persistence in nonautonomous delay differential Kolmogorov-type population models, Rocky Mountain J. Math., 24, 1--22


\bibitem{Lebedev} Andrei Lebedev, Mikhail Belousov. (1998) Functional Differential Equations.  CRC Press.


\bibitem{leung 1979} A. Leung (1979) Conditions for global stability concerning a prey-predator model with delay effects, SIAM J. Appl. Math., 36 , 281-286

\bibitem{leung and zhou 1988} Leung A W,Zhou (1988) Q.Global stability for a large class of Volterra-Lotka type integrodifferential population delay equation.  Nonlinear Analysis TMS, 12:495-505.






\bibitem{Li and Yorke 1975} T.Y Li and J.A. Yorke (1975) Period three implies chaos. Amer. Math. Monthly 82, 985--992.














\bibitem{Macdonald1989} N. MacDonald (1989) Biological Delay Systems: Linear Stability Theory. Cambridge Univ. Press: Cambridge.

\bibitem{Murray I:343} M.C. Mackey and J.G. Milton (1988) Dynamical diseases. Ann. N. Y. Acad. Sci. 504, 16--32.





\bibitem{mallet-paret and Nussbaum 1986} Mallet-Paret, J., Nussbaum, R.(1986) Global continuation and asymptotic behavior for periodic solutions of a differential delay equation. Ann. Mat. Pura. Appl. 145, 33-128






\bibitem{Miller 1976} R.S. Miller (1976) Pattern and process in competition.  Adv. Ecol. Res. 4, 1--74.







\bibitem{Murray 1984} J.D. Murray (1984) Asymptotic Analysis, Second Edition. Springer-Verlag: Berlin-heidelberg-New York.



\bibitem{Murray I} J.D. Murray. (2002) Mathematical Biology I: An Introduction, Third Edition.  Springer-Verlag: New York. 




















\bibitem{Post and Travis} W.M. Post and C.C. Travis (1981) Global stability in ecological models with continuous time delays. In: T. Herdman, H. Stech and S. Rankin (Eds.), Integral and Functional Differential Equations. Dekker: New York.









\bibitem{Shipanov1939} G.V. Schipanov (1939) Theory and methods of design of automatic controllers. Automat. Remote Control 1, 49--56.











\bibitem{Skubach97} Alexander L. Skubachevskii. (1997) Elliptic Functional Differential Equations and Applications. 
Birkhäuser.











\bibitem{Stepan1989} G. Stepan (1989) Retarded Dynamical Systems: Stability and Characteristic Functions. Longman Scientific and Technical: UK.

\bibitem{Stech 1979} H.W. Stech (1979) The Hopf-bifurcation: A stability result and application. J. Math. Anal. Appl. 71, 525--546.















\bibitem{Wang and Ma 1991} W. Wang and Z. Ma (1991) Harmless delays for uniform persistence. J. Math. Anal. Appl. 158, 256--268.








\bibitem{Woerz-Busekros} A. Woerz-Busekros (1978) Global stability in ecological systems with continuous time delay. SIAM J. Appl. Math. 35, 123--134.


\bibitem{Wright 1955} E.M. Wright (1955) A non-linear difference-differential equation. J. Reine Angew. Math. 494,66--87.







\bibitem{Yan 1987} J. Yan (1987) Oscillation of solutions of first order delay differential equations. Nonlin. Anal. 11, 1279--1287.






\end{thebibliography}
\end{document}